\pdfoutput=1
\documentclass{article}
\usepackage[margin=25mm]{geometry}
\usepackage{amsmath}
\usepackage{amsfonts}
\usepackage{amssymb}
\usepackage{graphicx}
\usepackage[utf8]{inputenc}
\usepackage[T2A]{fontenc}
\usepackage[ukrainian,russian,english]{babel}
\usepackage[pdftex]{hyperref}
\usepackage{xcolor}
\usepackage[all]{xy}
\newtheorem{lemma}{Lemma}
\newtheorem{theorem}{Theorem}

\title{Morse functions with four critical points on immersed 2-spheres}

\author{Svitlana Bilun, Bohdana Hladysh, 
Alexandr Prishlyak and Mariia Roman}
\begin{document}

\maketitle
\begin{abstract}
   We investigate topological properties of simple Morse functions with 4 critical points on immersed 2-spheres.  To classify such functions, dual graph of immersion and Reeb graphs is used.   We  have found all possible structures of the functions:6 structures with 4 critical points on one 1-strata  component, 7 structures with two points on the 1-strata and two points on the 2-strata, 7 structures with two 1-stratas and a three-connected 2-strata, three structures with two 1-stratas and without a three-connected 2-strata.
\end{abstract}
\section*{Introduction}
According to Reeb theorem, the sphere is a unique up to homeomorphism closed surface on which a (Morse) function with two critical points exists. The surface obtained as an image of a sphere under immersion have a larger number of points. We  show that there are at least four critical points on its.
At the same time, the surfaces we are considering do not have triple points.

%The main goal of the work is to describe the topological structures of such surfaces and functions with four critical points on them.

By Morse theory, an arbitrary Morse function on a manifold can be reduced to an optimal one by adding and reducing critical points. Therefore, it is important to know what structures are possible in optimal functions. Possible structures of optimal Morse functions on surfaces with a boundary are described in \cite{borodzik2016morse, hladysh2020deformations}

%A sphere immersed in a three-dimensional manifold has the natural structure of a stratified set: 1-strata are the components of the set of double points, and 2-strata are the components of their complements. A function $f: M \to R$ is called a \textit{Morse function} on a stratified set if every critical point belonging to the 2-strat is nondegenerate, and the critical points of the restriction of the function to the 1-strat are nondegenerate critical points as for the restriction to the 1-strat, as well as for limiting the function to the closure of the components of the intersection of 2-strats with the point \cite{Goresky1988}.

The use of graph theory to classify functions and dynamical systems is typical in two and three dimensions. 

%For Morse flows on surfaces, the first invariants were called distinguishing graphs \cite{Andronov1937, Leontovich1955, Smale1960, Peixoto1959, Peixoto1973, Palis1982, Palis1970, Fleitas1975}. After \cite{Oshemkov1998}, three-color graphs became popular \cite{prishlyak2020three}. Then these invariants were generalized to a wider class of dynamical systems, and other invariants were also constructed \cite{prislyak2017morse, Prishlyak2019, prishlyak2003topological, prishlyak2003sum,     Prishlyak2017, Prishlyak2022, Prishlyak2021,  Prishlyak2020, Kybalko2018,  Prishlyak2019,   Palis1968,  Giryk1996,   Bolsinov2004, Kadubovskyj2005, Poltavec1995}.

%For simple Morse functions on closed surfaces, the main invariant is the Reeb graph\cite{Reeb1946, Kronrod1950}. These graphs are also generalized to a wider class of functions \cite{prishlyak2001conjugacy, hladysh2019simple, hladysh2017topology, prishlyak2002topological1, prishlyak2002morse, prishlyak2000conjugacy,  prishlyak2007classification, lychak2009morse, bilun2002closed, 
%Sharko1993}.%Bilun13, 

%In dimension 3, diagrams consisting of a surface and the curves embedded in them are often used to classify flows and functions \cite{prishlyak1999equivalence,  prishlyak2002morse1,  Prishlyak2002, prishlyak2003regular, prishlyak2002topological, prishlyak2005complete, Prishlyak2007, Hatamian2020, BPP2022}. For polar Morse flows, they coincide with Heegaard diagrams.

%The main invariants of graphs and their embeddings in surfaces can be founded in \cite{prishlyak1997graphs, Harary69, HW68, GT87}.

The topological classifications of flows were obtained on closed 2-manifolds in \cite{bilun2023gradient, Kybalko2018, Oshemkov1998, Peixoto1973, prishlyak1997graphs, prishlyak2020three, akchurin2022three, prishlyak2022topological, prishlyak2017morse,  kkp2013, prishlyak2021flows,  prishlyak2020topology,   prishlyak2019optimal, prishlyak2022Boy}, 
 and on 2-manifolds with the boundary in
\cite{bilun2023discrete, bilun2023typical, loseva2016topology, loseva2022topological, prishlyak2017morse, prishlyak2022topological, prishlyak2003sum, prishlyak2003topological, prishlyak1997graphs, prishlyak2019optimal, stas2023structures}.
Complete topological invariants of Morse-Smale flows on 3-manifolds was constructed in \cite{prish1998vek,  prish2001top, Prishlyak2002beh2, prishlyak2002ms,  prishlyak2007complete, hatamian2020heegaard, bilun2022morse, bilun2022visualization}.

Morse flows are gradient flows of Morse functions in general position.  The flow determinate the topological structure of the function if one fixes the value of functions in singular points\cite{lychak2009morse, Smale1961}. Therefore, Morse--Smale flows classification is related to the classification of the Morse functions.

Topological invariants of functions on oriented 2-maniofolds were constructed in \cite{Kronrod1950} and \cite{Reeb1946} and  in \cite{lychak2009morse} for  non-orientable two-dimensional manifolds, in   \cite{Bolsinov2004, hladysh2017topology, hladysh2019simple, prishlyak2012topological} for manifolds with boundary, in \cite{prishlyak2002morse} for non-compact manifolds. 

Topological invariants of smooth function on closed 2-manifolds was also investigated in \cite{bilun2023morseRP2, hladysh2019simple, hladysh2017topology,  prishlyak2002morse, prishlyak2000conjugacy,  prishlyak2007classification, lychak2009morse, prishlyak2002ms, prish2015top, prish1998sopr,  bilun2002closed,  Sharko1993}, on 2-manifolds with the boundary in \cite{hladysh2016functions, hladysh2019simple, hladysh2020deformations} and on closed 3- and 4-manifolds in  \cite{prishlyak1999equivalence, prishlyak2001conjugacy}.

For the first acquaintance with the theory of functions and flows on manifolds, we recommend the papers \cite{prishlyak2012topological, prish2002theory, prish2004difgeom, prish2006osnovy, prish2015top}.
%Catastrophe saddle-saddle  $x^3 - t x - x y^2$, %{x, -1, 1.5}, {y, -2.5, 2.5}

%Bifurcation saddle-saddle $\{ 3 x^2 - t - y^2, -2 x  y \}$

\textbf{Statement of the problem.}
The purpose of this paper is to describe all possible topological structures of simple Morse functions with four critical points on stratified sets, which is the immersion of a two-dimensional sphere into a three-dimensional oriented manifold without triple points. 
 
\textbf{Basic Concepts.}
A stratified set is a space divided into subsets (strata), each of which is a manifold, and each point has a neighbourhood that intersects with a finite number of strata. When a compact two-dimensional manifold is immersed in a three-dimensional manifold, there is a natural stratification: strata of dimension 0 (0-strata) are the components of sets of triple points (the pre-image contains three different points), 1-strata are the components of the set of double points (the pre-image is 2-points), 2 - strata are components of a set of points, the pre-image of which is a single point.

Let $M$ be a two-dimensional stratified set without 0-strata. A function $f: M \to \mathbb{R}$ is called \textit{a Morse function} if every critical point belonging to the 2-strata is non-degenerate, and the critical points of the restriction of the function to the 1-strata are non-degenerate critical points as for the restriction to the 1-strata, as well as for limiting the function to the closure of the components of the intersection of 2-stratas with the neighbourhood of a point \cite{GM88}.

The Morse function is called \textit{simple} if it takes different values at different critical points. In other words: each critical value corresponds to one critical point.

%The Morse function $f: M \to \mathbb{R}$ is called optimal if there is no Morse function on $M$ with a smaller number of critical points than in $f$.
\textbf{Paper structure.}
In the first section, to describe the structure of the immersion of the sphere, we introduce a graph that is a dual graph to the pre-image of 1-strata on the sphere, the edges of which are divided into pairs (colouring) depending on the coincidence of the images during immersion. Since each 1-strata contains at least two critical points (points of minimum and maximum for restriction of the function per the 1-strata), the total number of 1-stratas in a function with four critical points does not exceed two. Therefore, all possible stratifications with one and two 1-strata are considered.

In the second section, we use Reeb graphs for the topological classification of function restrictions on 2-strata closures.

In the third section, we use the results of the first two sections to describe all possible structures of Morse functions with four critical points on immersed spheres in the case when the set of double points is connected, and in the fourth section, when it has two components.

\section{Dual  graph of the immersion}

If $g: S^2 \to N^3$ is an immersion of the sphere, then the stratification $M=g(S^2)$ defines the stratification of the sphere, where the stratas are the components of the pre-images of the stratas. If the natural stratification $M$ had a single 1-stratum, then the corresponding sphere stratification have two 1-stratas. If the natural stratification $M$ had two 1-strats, then there are four 1-strats on the sphere. We use dual graphs to describe their possible configurations.

\textit{A dual graph} of stratification of the sphere is called a graph whose vertices correspond to 2-strata and edges to 1-strata. At the same time, two vertices 
are connected by an edge if the corresponding circle (1-strata) lies in the boundary of each of the 2-stratas corresponding to these two vertices.
%By analogy with \cite{Palm2012} is the case
\begin{lemma}
A dual stratification graph on a sphere is a tree. Two stratifications of a sphere are homeomorphic if and only if their dual graphs are isomorphic.
\end{lemma}

The proof of the lemma is analogous to the proof of the corresponding fact in the classification of functions on three-dimensional manifolds in \cite{prishlyak2012topological}.

\begin{figure}[ht!]
\center{\includegraphics[width=0.75\linewidth]{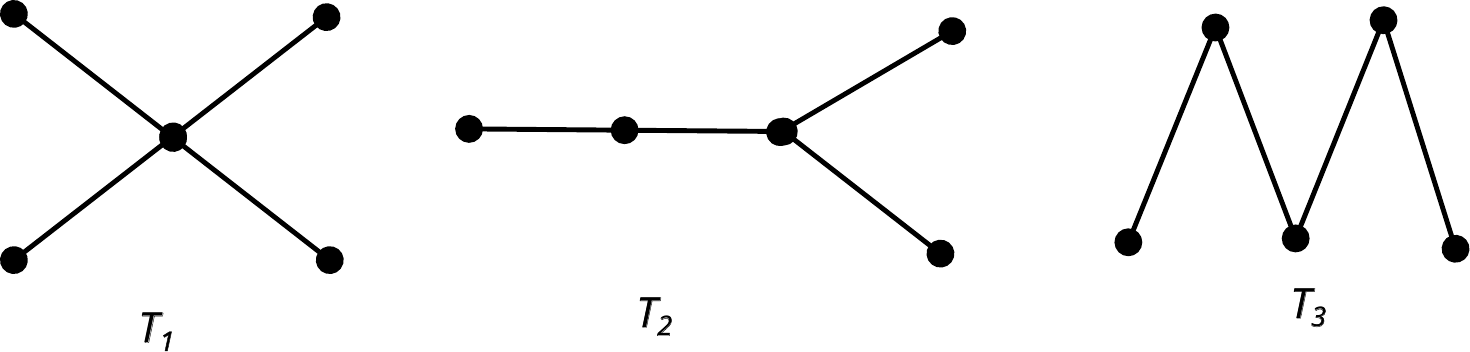}
}
\caption{trees with four edges
}
\label{trees}
\end{figure}

If the stratification on the sphere has two 1-stratas, then the dual graph is a tree  and has two edges. Since there is only one tree with two edges, all stratifications of a sphere with two edges are homeomorphic.

If the sphere stratification has four 1-strata, then the dual graph has 4 edges. There are three trees with four edges (see Fig. \ref{trees}):

1) $T_1$ has one vertex of degree 4, and the other vertices have degree 1 (in other words, all edges have a common vertex);

2) $T_2$ has one vertex of degree 3, one vertex of degree 2, and three vertices of degree 1;

3) $T_3$ has three vertices of degree 2 and two vertices of degree 1 (chain).

\begin{figure}[ht!]
\center{\includegraphics[width=0.85\linewidth]{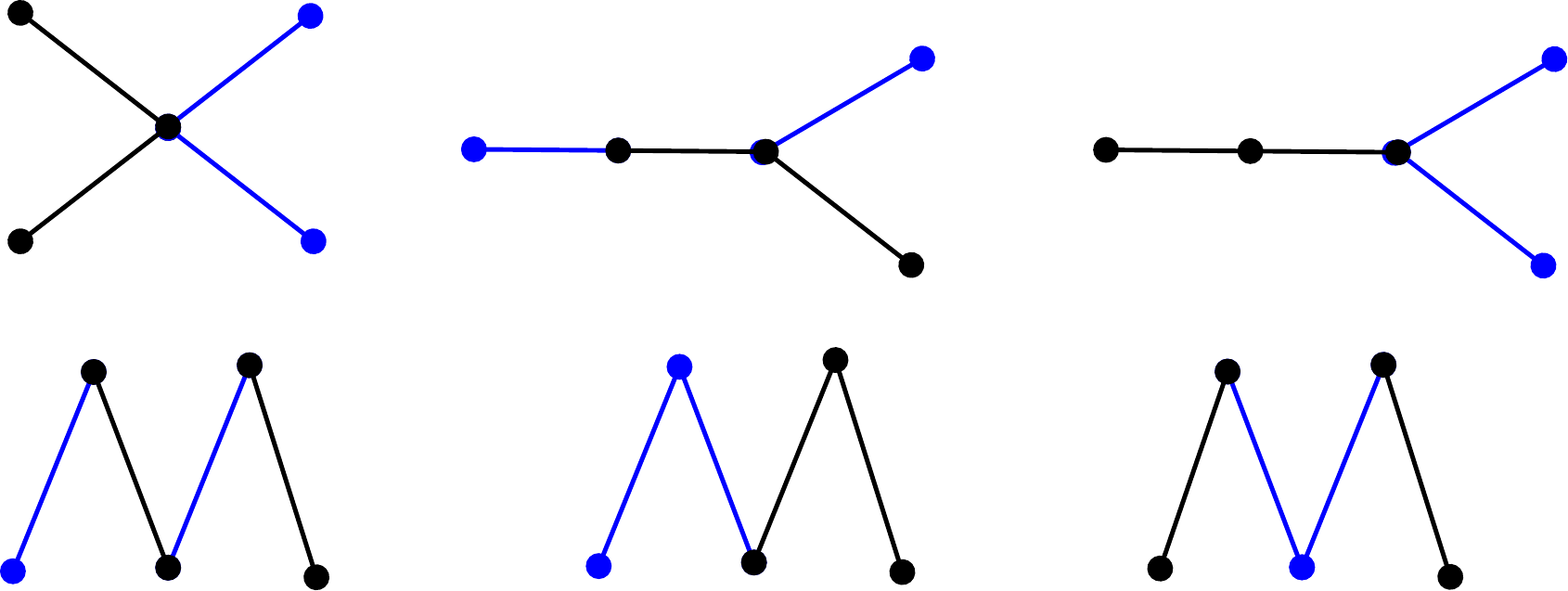}
}
\caption{coloring of trees with four edges
}
\label{col}
\end{figure}

Let's color the edges of the tree in two colors so that two edges have one color, if the corresponding stratas are displayed in one strata when the sphere is immersed. Then, with accuracy up to isomorphism, the graph $T_1$ has a single coloring, $T_2$ has two colorings, and %It is necessary to describe why the second coloring of the second graph is not possible
graph $T_3$ has three colorings, see Fig. \ref{col}

%A) the central vertex divides the graph into subgraphs colored in the same color;

%B) the colors of the edges in the graph alternate.

Let's call the vertices of degree 1 marginal, and the remaining vertices internal. An edge is marginal if one of its ends has degree 1, otherwise the edge is internal. We also call the corresponding strata internal and marginal.

\section{Reeb graph of a simple Morse function on an oriented surface with embedded  circles}

Let the stratified set be a compact oriented surface, the 1-strata are simply closed curves and surface boundary components, and the Morse function has exactly two critical points on each strata. If we glue together the points of each component of the boundary with the same values, then we get a closed surface with segments embedded in it.

The Reeb graph is an  graph, which is a factor space of the surface according to the equivalence relation in which points are equivalent if they belong to the same component of the level of the function. The orientation of the edges is determined by the direction of growth of the function.

\begin{figure}[ht!]
\center{\includegraphics[width=0.65\linewidth]{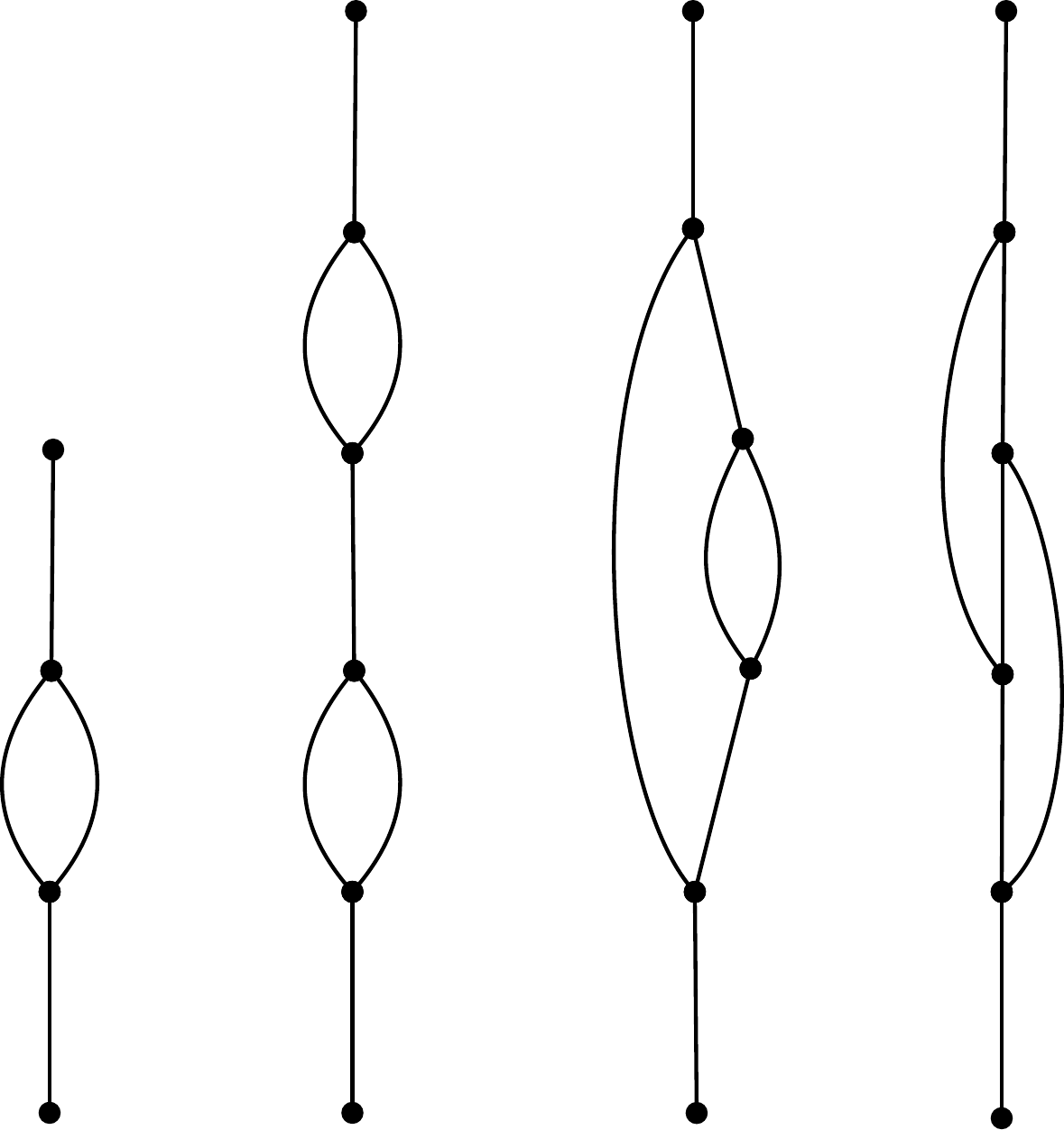}
}
\center{
1) \ \ \ \ \ \ \ \ \ \ \ \ \ \ \ \ \ \ \ \ \ \ \ \
2) \ \ \ \ \ \ \ \ \ \ \ \ \ \ \ \ \ \ \ \ \ \ \ \ \ \
3)\ \ \ \ \ \ \ \ \ \ \ \ \ \ \ \ \ \ \ \ \ \ \ \ \
4)}
\caption{ Rib graphs of optimal Morse functions on 1) torus and 2) -- 4) double torus
}
\label{reebs}
\end{figure}

If there are no 1-strata, then the Reeb graph is a complete topological invariant of a simple Morse function. Each 1-strata is divided by critical points into two paths such that the value of the function increases when moving along each of them. These paths generate paths on the Reeb graph. We number the paths so that 1 and 2, 3 and 4, 5 and 6, etc., correspond to one 1-stratum, and the remaining paths that are not divided into pairs correspond to segments obtained after gluing the points of the boundary components with the same function values.

Fix the orientation of the surface. At each level component, the orientation is induced from the orientation of the set of points having a smaller function value than this level. A critical level component containing a critical point is split by the critical point into segments. The direction of movement based on the orientation of this segment determines the %license 
order on the set of points of intersection of it with paths. Therefore, the set of paths, the Reeb graph, passing through the vertex of valence 3 is divided into two ordered subsets. These orders induce cyclic orders with respect to edges identical to them. Note that if both ends of an edge have a valency of 3, then the two induced cyclic orders must coincide.

\textit{A distinguishing graph} is a Reeb graph with paths selected on it, some of which are split into pairs, and for each vertex of valence 3, a partition of the set of paths passing through it into two ordered subsets.

\begin{theorem} Let $M$ be a compact oriented surface on which the structure of a stratified set without 0-stratas is given, $f,g: M \to R$ are Morse functions on this stratified set, which have two critical points on to each 1-stratum. The functions $f$ and $g$ are fibre equivalent if and only if their distinguishing graphs are equivalent.
\end{theorem}

\textbf{Proof} of the theorem is similar to the proof of the corresponding theorems in 
\cite{Kronrod1950, Reeb1946, Bolsinov2004}. 

\textbf{Properties of the distinguishing graph.}
Paths belonging to one of the two subgroups at the vertex pass through one of the two adjacent edges included
into this vertex (or leave it).

If a vertex is the start or end of two paths (paths from the same pair), then these paths belong to different subsets in this vertex.

These properties follow directly from the definition.

\section{Functions with a single 1-stratum on an immersed sphere}

Let us first consider the case when there are two critical points on the 1-stratum.

The apex of step 1 corresponds to a single-connected region (2-disc). The optimal Morse function on it has two critical points - a minimum and
maximum (which are on the 2-disk boundary). Therefore, with accuracy up to topological equivalence, there is a single optimal function on such 2-strata.

If the 2-stratum has two components of the border, then in the image when immersed, they are displayed in one 1-stratum, so the function takes the same values on them. If you glue them together (take the closure of the image of this stratum during immersion), then you will get a torus. The optimal Morse function on a torus has 4 critical points: the minimum point $p_0$, the first saddle $p_1$, the second saddle $p_2$ and the maximum point $p_3$. Correspondent Reeb graphs is shown in fig. \ref{reebs} 1). The two critical points on the 1-strata will be two of these four critical points. The following pairs of points are possible:
\begin{center}
1) $(p_0,p_1)$, 2) $(p_0,p_2)$, 3) $(p_0,p_3)$, 4) $(p_1,p_2)$, 5) $(p_1,p_3)$, 6) $(p_2,p_3)$. 
\end{center}
For the first pair, both paths going from $p_0$ to $p_1$ coincide on the Reeb graph and form a meridian on the torus. Therefore, there is a single Morse function structure in this case. In the second case, we have two different paths on the Reeb graph that form a parallel on the torus. In the third case, two pairs of paths are possible: a) they coincide on the Reeb graph and form a meridian on the torus, b) different on the Reeb graph are parallel on the torus. In the fourth case, two different paths form a parallel on the torus. The fifth and sixth cases are analogous to the second and first cases for the inverse Morse function.

$ $

Now consider the case when all 4 critical points lie on the 1-stratum. For each of the 2-strata, homeomorphic to 2-disks, the Reeb graph has the form of the letter Y or its inverted letter. There are thus two structures of Morse functions on each of them, but taking into account symmetry, there are three structures on a pair of these 2-discs. On the torus, the 1-strata passes through all four critical points. If you move along the 1-stratum, starting from the minimum point to the critical point with the smallest value, then only the following sequence of critical points is possible: $p_0, p_2, p_1, p_3$. At the same time, on the Reeb graph, the image of this cycle either covers the entire graph or leaves one of the two edges between $p_1$ and $p_2$ uncovered. Therefore, we have two structures on the torus with a %nested 
loop. In general, on the immersed sphere there are $2 \times 3=6$ structures with four critical points on one 1-strata.

Summarizing all of the above, we have the following statement:

\begin{theorem}
There are 13 topologically non-equivalent optimal Morse functions on an immersed sphere with a connected set of double points.
\end{theorem}

\section{Functions with two 1-strata on an immersed sphere}

\begin{figure}[ht!]
\center{\includegraphics[width=0.95\linewidth]{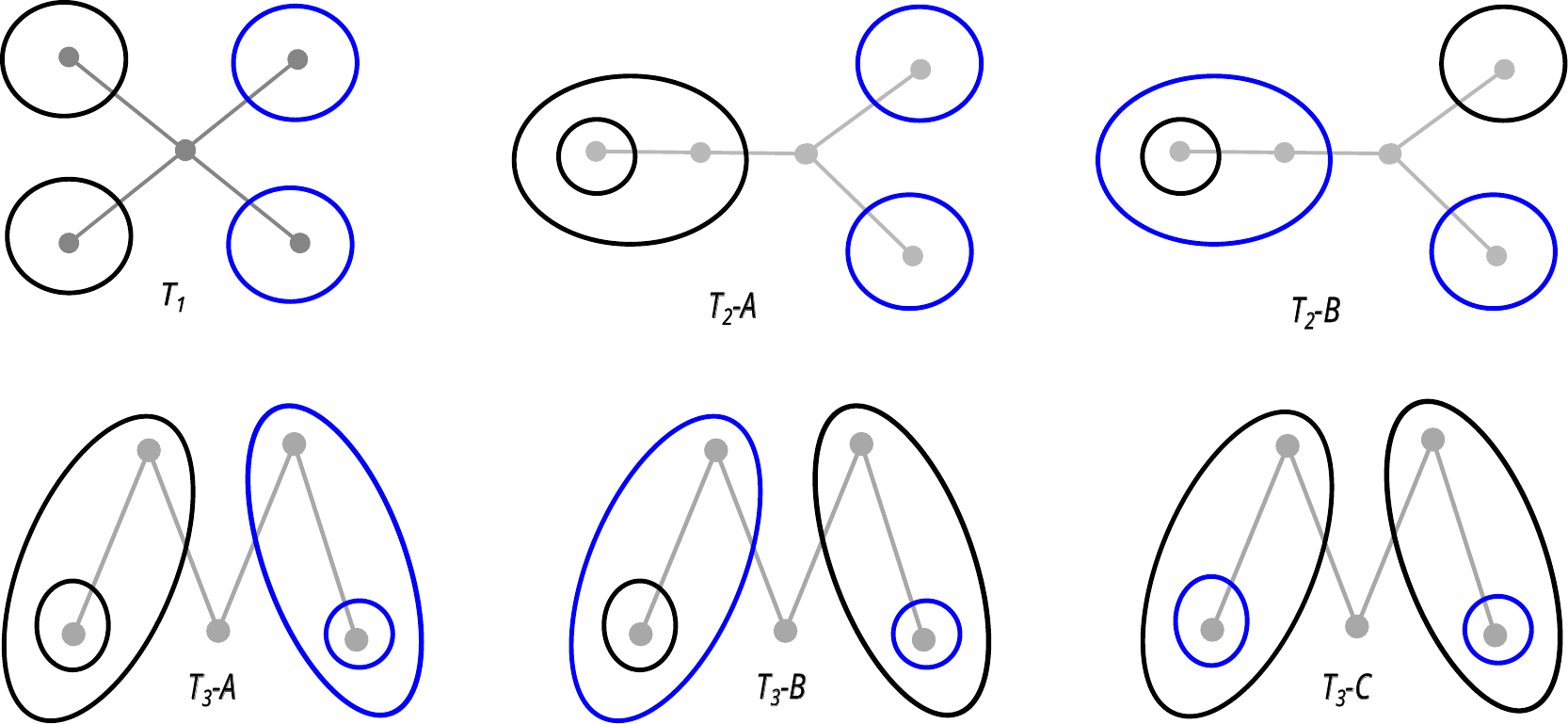}
}
\caption{stratification of submerged spheres
}
\label{strat}
\end{figure}

In fig. \ref{strat}, all possible stratifications of immersed spheres with two components of the set of double points are shown. They correspond to the colouring of the graphs described in fig. \ref{col}. Let's consider the Morse functions on each of them in details.

\textbf{1. Graph $T_1$.}
The closure of the image of the inner 2-strata is an oriented surface of genus 2. The optimal Morse function on such a surface has 6 critical points, ordered by the growth of the function on them. Correspondent Reeb graphs are shown in fig. \ref{reebs} 2)--4). Therefore, there are no Morse functions with four critical points on such a stratified set.

\textbf{2. Graph $T_2 -A$.}
The closure of the interior regions is a torus and a torus with a hole. The torus corresponds to the vertex of valence 2 on the dual graph. On this torus, in addition to the two critical points on the 1-strata, there are two more internal critical points. That is, the total number of critical points is at least six.

Therefore, there are no Morse functions with four critical points on such a stratified set.

\textbf{3. Graph $T_2 -B$.}
On two-dimensional disks (corresponding to vertices of valence 1 on the dual graph), as well as on the cylinder (corresponding to vertex of valence 2), the Morse function on the boundary completely determines the structure of the Morse function inside, since there is always a single such continuation without internal critical points. Therefore, the structure of the Morse function on this stratified set is completely determined by the structure on the stratum corresponding to the vertex of valence 3 on the dual graph. For the blue cycle, as well as in the case of one 1-strata, 7 cases are possible. In all cases, except for case 3 a), there exists a symmetry of the Reeb graph that preserves the pair of paths. Thus, in cases with symmetry, the two remaining vertices can be connected by a single path, with accuracy up to homeomorphism. In case 3 a), there are two such possibilities for choosing a path (whether or not the union of all paths covers the Reeb graph). So there are eight structures of Morse functions with four critical points in this case.

\textbf{4. Graph $T_3 -A$.}
In this case, the closure of two internal 2-stratas is a torus, and another 2-strata is a cylinder.
On each of the torus, the optimal function has 4 critical points. Therefore, the total number of an arbitrary Morse function on this surface is not less than 8.

Therefore, there are no Morse functions with four critical points on such a stratified set.

\textbf{5. Graph $T_3-B$.}
In this case, the closure of each internal 2-strata is a cylinder. On each of them, there is a single structure of the optimal Morse function. Critical values on 1-strata may alternate, or both values on one stratum may be less than critical values on the other, or one stratum may have a minimum and maximum point, and the others has two average critical points.

So, there are 3 structures in this case.

\textbf{6. Graph $T_3 -C$.}
In this case, the closures of two internal 2-stratas are cylinders, and another 2-strata is a torus.
On the torus, the optimal function has 4 critical points. Therefore, the total number of an arbitrary Morse function on this surface is at least 6.

Therefore, there are no Morse functions with four critical points on such a stratified set.

$ $

Summarizing all of the above, we have the following statement:

\begin{theorem}
On the immersed spheres with two components of the set of dual points, there are 11 topologically non-equivalent optimal Morse functions.
\end{theorem}

\section*{Conclusion}
We have described all possible structures of sphere immersions in three-dimensional manifolds with a connected set of double points and two components of the set of double points. On each of these spaces, all possible structures of simple Morse functions with four critical points are found: 13 structures for the connected set of dual points and 11 structures for the two connected components. To classify these structures, we use Reeb graphs for functions on the closure of each 2-strata.

We hope that the obtained results can be generalized to functions with a larger number of critical points on immersed spheres and other surfaces. Since the constructed structures are discrete in nature, it also is interesting to describe all possible structures of discrete Morse functions on these and other stratified two-dimensional sets.

%\bibliographystyle{plain}
%\bibliography{prishd}

\textsc{Taras Shevchenko National University of Kyiv}

Svitlana Bilun  \ \ \ \ \ \ \ \ \ \ \ \
\textit{Email address:} \text{ bilun@knu.ua}   \ \ \ \ \ \ \ \ \ \
\textit{ Orcid ID:}  \text{0000-0003-2925-5392}

Bohdana Hladysh \ \ \ \ \ \ \ 
\textit{Email:} \text{ bohdanahladysh@gmail.com,}  \ 
\textit{ Orcid ID:} \text{0000-0001-8935-1453}

Alexandr Prishlyak \ \ \ \ \ \ 
\textit{Email address:} \text{ prishlyak@knu.ua} \ \ \ \  \
\textit{ Orcid ID:} \text{0000-0002-7164-807X}

Mariia Roman \ \  \ \ \ \ \ \ \ \ \ \  \textit{Email address:} \text{3848855@gmail.com }    \ \
\textit{ Orcid ID:} \text{0009-0000-6795-8991}

\end{document}